\documentclass[11pt]{amsart} 
\usepackage{amsmath,amssymb}
\usepackage[dvips]{graphicx}  
\usepackage{psfrag}
\def\COMMENT#1{}
\def\TASK#1{}   

\def\noproof{{\unskip\nobreak\hfill\penalty50\hskip2em\hbox{}\nobreak\hfill%
        $\square$\parfillskip=0pt\finalhyphendemerits=0\par}\goodbreak}
\def\endproof{\noproof\bigskip}
\newdimen\margin   
\def\textno#1&#2\par{%
    \margin=\hsize
    \advance\margin by -4\parindent
           \setbox1=\hbox{\sl#1}%
    \ifdim\wd1 < \margin
       $$\box1\eqno#2$$%
    \else
       \bigbreak
       \hbox to \hsize{\indent$\vcenter{\advance\hsize by -3\parindent
       \sl\noindent#1}\hfil#2$}%
       \bigbreak
    \fi}
\def\proof{\removelastskip\penalty55\medskip\noindent{\bf Proof. }}

\def\cP{\mathcal{P}}
\def\E{\mathcal{E}}
\def\F{\mathcal{F}}
\def\eps{\varepsilon}

\newtheorem{firstthm}{Proposition}
\newtheorem{thm}[firstthm]{Theorem}

\newtheorem{lemma}[firstthm]{Lemma}

\newtheorem{conj}[firstthm]{Conjecture}
\newtheorem{claim}{Claim}   
\newtheorem{fact}[firstthm]{Fact}

\begin{document}
\title{Matchings in $3$-uniform hypergraphs}
\author{Daniela K\"uhn, Deryk Osthus and Andrew Treglown}
\date{\today} 

\begin{abstract} 
We determine the minimum vertex degree that ensures a perfect matching in a $3$-uniform hypergraph. More precisely,
suppose that $H$ is a sufficiently large $3$-uniform hypergraph whose order $n$ is divisible by $3$. If
the minimum vertex degree of $H$ is greater than $\binom{n-1}{2}-\binom{2n/3}{2}$,
then $H$ contains a perfect matching. This bound is tight and answers a question of H\`an, Person and Schacht.
More generally, we show that $H$ contains a matching of size $d\le n/3$ if its minimum vertex degree is greater than
$\binom{n-1}{2}-\binom{n-d}{2}$, which is also best possible.
This extends a result of Bollob\'as, Daykin and Erd\H{o}s. 
\end{abstract}

\maketitle
\section{Introduction}\label{sec1}
A \emph{perfect matching} in a hypergraph $H$ is a collection of vertex-disjoint edges of $H$ which cover the vertex set $V(H)$ of $H$. A theorem of 
Tutte~\cite{tutte} gives a characterisation of all those graphs which contain a perfect matching. On the other hand, the decision problem whether
an $r$-uniform hypergraph contains a perfect matching is NP-complete for $r\geq 3$. (See, for example,~\cite{ruc2} for complexity results in the
area.) It is natural therefore to seek simple sufficient conditions, such as minimum degree conditions,
that ensure a perfect matching in an $r$-uniform hypergraph.
This has turned out to be a difficult question: despite considerable attention, the full solution remains elusive.
But the partial results obtained so far have already involved the development of new techniques and uncovered interesting connections to other problems. 

Given an $r$-uniform hypergraph $H$ and distinct vertices $v_1, \dots , v_{\ell} \in V(H)$ (where $1 \leq \ell \leq r-1$) we define 
$d_H (v_1, \dots , v_{\ell})$ to be the number of edges containing each of $v_1, \dots , v_{\ell}$. The \emph{minimum $\ell$-degree $\delta _{\ell} 
(H)$} of $H$ is the minimum of $d_H (v_1, \dots , v_{\ell})$ over all $\ell$-element sets of vertices in $H$. Of these parameters the two most natural
to consider are the \emph{minimum vertex degree} $\delta _1 (H)$ and the \emph{minimum collective degree} or \emph{minimum codegree} $\delta _{r-1}
(H)$. R\"odl, Ruci\'nski and Szemer\'edi~\cite{rrs} determined the minimum codegree that ensures  a perfect matching in
an $r$-uniform hypergraph. 
This improved bounds given in~\cite{ko1, rrs2}. An $r$-partite version was proved by
Aharoni, Georgakopoulos and Spr\"ussel~\cite{ags}.

Much less is known about minimum vertex degree conditions for perfect matchings in $r$-uniform hypergraphs $H$. 
H\`an, Person and Schacht~\cite{hps} showed that the threshold in the case when $r=3$ is $(1+o(1))\frac{5}{9}\binom{|H|}{2}$. 
(Here, $|H|$ denotes the number of vertices in $H$.)
This improved an earlier bound given by Daykin and H\"aggkvist~\cite{dayhag}. 
In this paper we determine the threshold exactly, which answers a question from~\cite{hps}.

\begin{thm}\label{mainthm} There exists an $n_0 \in \mathbb N$ such that the following holds.
Suppose that $H$ is a $3$-uniform hypergraph whose order $n \geq n_0$ is divisible by $3$. If
$$\delta _1 (H) >  \binom{n-1}{2}-\binom{2n/3}{2} $$
then $H$ has a perfect matching.
\end{thm}
Independently, Khan~\cite{khan1} has given a proof of Theorem~\ref{mainthm} using different arguments.
The following example shows that the result is best possible: let $H^*$ be the $3$-uniform hypergraph whose vertex set is
partitioned into two vertex classes $V$ and $W$ of sizes $2n/3 +1$ and $n/3-1$ respectively and whose edge set consists precisely
of all those edges with at least one endpoint in $W$. Then $H^*$ does not have a perfect matching and
$\delta _1 (H) =  \binom{n-1}{2}-\binom{2n/3}{2}$.

The example generalises in the obvious way to $r$-uniform hypergraphs. This leads to the following 
conjecture, which is implicit in several earlier papers (see e.g.~\cite{hps,survey}).
Partial results were proved by H\`an, Person and Schacht~\cite{hps} as well as Markstr\"om and Ruci\'nski~\cite{mark}.
\begin{conj}\label{conjj}
For each integer $r\ge 3$ there exists an integer $n_0=n_0(r)$ such 
that the following holds.
Suppose that $H$ is an $r$-uniform hypergraph whose order $n \geq n_0$ is divisible by $r$. If
$$\delta _1 (H) >  \binom{n-1}{r-1}-\binom{(r-1)n/r}{r-1},$$
then $H$ has a perfect matching. 
\end{conj}
Recently, Khan~\cite{khan2} proved Conjecture~\ref{conjj} in the case when $r=4$.
It is also natural to ask about the minimum (vertex) degree which guarantees a matching of given size $d$.
Bollob\'as, Daykin and Erd\H{o}s~\cite{bde} solved this problem for the case when $d$ is small compared to the 
order of~$H$. We state the $3$-uniform case of their result here. The above hypergraph $H^*$ with $W$ of size $d-1$
shows that the minimum degree bound is best possible. 
\begin{thm}[Bollob\'as, Daykin and Erd\H{o}s~\cite{bde}]\label{bde}
Let $d \in \mathbb N$. If $H$ is a $3$-uniform hypergraph on $n>54(d+1)$ vertices and
$$\delta _1 (H) >\binom{n-1}{2} - \binom{n-d}{2}$$
then $H$ contains a matching of size at least $d$.
\end{thm}
Here we extend this result to the entire range of $d$. Note that Theorem~\ref{generalmatch} generalises Theorem~\ref{mainthm},
so it suffices to prove Theorem~\ref{generalmatch}.
\begin{thm}\label{generalmatch}
There exists an $n_0 \in \mathbb N$ such that the following holds.
Suppose that $H$ is a $3$-uniform hypergraph on $n\ge n_0$ vertices, that $n/3\ge d\in\mathbb{N}$ and that
$$\delta _1 (H) >\binom{n-1}{2} - \binom{n-d}{2}.$$
Then $H$ contains a matching of size at least $d$.
\end{thm}

It would be interesting to obtain analogous results (i.e. minimum degree conditions which guarantee a matching
of size $d$) for $r$-uniform hypergraphs and for $r$-partite hypergraphs. Some bounds are given in~\cite{dayhag}. 
Further, a $3$-partite version of Theorem~\ref{mainthm} was recently proved by Lo and  Markstr\"om~\cite{lomark}.

Treglown and Zhao~\cite{zhao, zhao2} determined the minimum $\ell$-degree that ensures a perfect matching in an
$r$-uniform hypergraph when $r/2 \leq \ell \leq r-1$. (Independently, Czygrinow and Kamat~\cite{czy} dealt with the
case when $r=4$ and $\ell=2$.) Prior to this, 
Pikhurko~\cite{pik} gave an asymptotically exact result.
The situation for $\ell$-degrees where $1 < \ell < r/2$ is still open.
In~\cite{hps}, H\`an, Person and Schacht provided conditions on $\delta _{\ell} (H)$ that ensure a perfect matching in the case
when $\ell< r/2$. These bounds were subsequently lowered by Markstr\"om and Ruci\'nski~\cite{mark}.
 Alon, Frankl, Huang, R\"odl, Ruci\'nski and Sudakov~\cite{afh} discovered 
a connection between the minimum $\ell$-degree
that forces a perfect matching in an $r$-uniform hypergraph and the  minimum $\ell$-degree
that forces a \emph{perfect fractional matching}. As a consequence of this result they determined, asymptotically,
the minimum $\ell$-degree that ensures a perfect matching in an $r$-uniform hypergraph for the following values of
$(r, \ell)$: $(4, 1)$, $(5, 1)$, $(5, 2)$, $(6, 2)$ and  $(7, 3)$. 
See~\cite{rrsurvey} for further results concerning perfect matchings in hypergraphs.

\section{Notation}\label{sec2}

Given a hypergraph $H$ and  subsets $V_1,V_2,V_3$ of its vertex set $V(H)$, we say that an edge $v_1v_2v_3$ is of \emph{type}
$V_1V_2V_3$ if $v_1 \in V_1$, $v_2 \in V_2$ and  $v_3 \in V_3$.%
   \COMMENT{Also had: We also write $e_H(V,W)$ for the number of all
those edges which are of type $VVW$ or $VWW$. Check that this is not needed.}

Let $d\le n/3$ and let $V,W$ be a partition of a set of $n$ vertices such that $|W|=d$.
Define $H_{n,d}(V,W)$ to be the hypergraph with vertex set $V \cup W$ consisting of all those edges which
have type $VVW$ or $VWW$. Thus $H_{n,d}(V,W)$ has a matching of size~$d$,
$$\delta_1(H_{n,d}(V,W))= \binom{n-1}{2} - \binom{n-d-1}{2}$$
and $H_{n,d}(V,W)$ is very close to the extremal hypergraph which shows that
the degree condition in Theorem~\ref{generalmatch} is best possible.
$V$ and $W$ are the \emph{vertex classes of $H_{n,d}(V,W)$}.

Given $\eps >0$, a $3$-uniform hypergraph $H$ on $n$ vertices and a partition $V,W$ of $V(H)$ with $|W|=d$, we say that $H$ is \emph{$\eps$-close
to $H_{n,d}(V,W)$} if  
$$|E(H_{n,d}(V,W)) \setminus E(H)| \le \eps n^3.$$
In this case we also call $V$ and $W$ \emph{vertex classes of $H$}.
(So $H$ does not have unique vertex classes.) We say that $H$ is \emph{$\eps$-close
to $H_{n,d}$} if there is a partition $V,W$ of $V(H)$ such that $|W|=d$ and $H$ is $\eps$-close
to $H_{n,d}(V,W)$.

Given a vertex $v$ of a $3$-uniform hypergraph $H$, we write $N_H(v)$ for the \emph{neighbourhood of $v$},
i.e.~the set of all those (unordered) tuples of vertices which form an edge together with~$v$.
Given two disjoint sets $A,B\subseteq V(H)$, we define the \emph{link graph $L_v(A,B)$
of $v$ with respect to $A,B$} to be the bipartite graph whose vertex classes are $A$ and $B$ and in which
$a\in A$ is joined to $b\in B$ if and only if $ab\in N_H(v)$. 
Similarly, given a set $A\subseteq V(H)$, we define the \emph{link graph $L_v(A)$
of $v$ with respect to $A$} to be the graph whose vertex set is $A$ and in which
$a,a'\in A$ are joined if and only if $aa'\in N_H(v)$. 
Also, given disjoint sets $A,B,C,D,E\subseteq V(H)$, we write
$L_v (ABCD)$ for $L_v (A,B) \cup L_v (B,C) \cup L_v (C,D)$. We define $L_v (ABCDE)$ similarly.
If $M$ is a matching in $H$ and $E,F$ are two edges in $M$ with $v\notin E,F$, we write $L_v(EF)$ for
$L_v(V(E),V(F))$.  If $E_1,\dots,E_5$ are matching edges avoiding $v$, we define $L_v(E_1\dots E_4)$ and $L_v(E_1\dots E_5)$ similarly.
If $e=uw$ is an edge in the link graph of $v$, then we write $ve$ for the edge $vuw$ of $H$.
A matching in $H$ of size $d$ is called a \emph{$d$-matching}.

Given a set $M$ and $k\ge 2$, we write $\binom{M}{k}$ for the set of all $k$-element subsets of $M$%
\COMMENT{and $(M)^k$ for the set of all ordered $k$-tuples of elements of $M$ - not needed}. Given sets $M$ and $M'$, we write
$MM'$ for the set of all pairs $mm'$ with $m\in M$ and $m'\in M'$.

Given two graphs $G$ and $G'$, we write $G\cong G'$ if they are isomorphic.
A bipartite graph is called \emph{balanced} if its vertex classes have equal size.
By a \emph{directed graph} we mean a graph whose edges are directed, but we only allow at most two edges between
any pair of vertices: at most one edge in each direction. We write $vw$ for the edge directed from $v$ to $w$. 
Given disjoint vertex sets $V$ and $W$ of a directed graph, we write $e(V,W)$ for the number of all those edges
which are directed from some vertex in $V$ to some vertex in $W$. A directed graph $G$ is an \emph{oriented graph} if
it has at most one edge between any pair of vertices (i.e.~if $G$ has no directed cycle of length~2).

We will often write $0<a_1 \ll a_2 \ll a_3$ to mean that we can choose the constants
$a_1,a_2,a_3$ from right to left. More
precisely, there are increasing functions $f$ and $g$ such that, given
$a_3$, whenever we choose some $a_2 \leq f(a_3)$ and $a_1 \leq g(a_2)$, all
calculations needed in our proof are valid. 
Hierarchies with more constants are defined in the obvious way.

\section{Preliminaries and outline of proof}

Our approach towards Theorem~\ref{generalmatch} follows the so-called \emph{stability approach}: 
we prove an approximate version of the desired result which states that the minimum degree condition implies
that either (i) $H$ contains a $d$-matching or (ii) $H$ is `close' to the extremal hypergraph. The latter implies that $H$ is `close'
to the hypergraph $H_{n,d}$ defined in the previous section. This extremal situation (ii) is then dealt with separately. 
We do this in Section~\ref{extremal}, where we prove Lemma~\ref{extremallemma}.
The proof of Lemma~\ref{extremallemma} makes use of Theorem~\ref{bde}.

The non-extremal case is proved in Section~\ref{mainsec}.
As mentioned earlier, an approximate version of Theorem~\ref{mainthm} was proved in~\cite{hps}.
However, we need to proceed somewhat differently as the argument in~\cite{hps} fails to guarantee
the `closeness' of $H$ to the extremal hypergraph in case (ii).
(But we do use the same general approach and a number of ideas from~\cite{hps}.)

We begin by considering a matching $M$ of maximum size and suppose that $|M|<d$.
We then carry out a sequence of steps, where in each step we show that we can either find a larger matching
(and thus obtain a contradiction), or show that $H$ is successively `closer' to $H_{n,d}$.
Amongst others, the following fact from~\cite{hps} will be used to achieve this
(see Figure~1 for the definitions of $B_{033},B_{023}, B_{113}$).
\begin{fact}\label{balanced}
Let $B$ be a balanced bipartite graph on $6$ vertices. 
\begin{itemize} 
\item If $e(B)\ge 7$ then $B$ contains a perfect matching.
\item If $e(B)=6$ then either $B$ contains a perfect matching or $B\cong B_{033}$.
\item If $e(B)=5$ then either $B$ contains a perfect matching or $B\cong B_{023}, B_{113}$.
\end{itemize}
\end{fact}
\begin{figure}[htb!]
\begin{center}\footnotesize
\includegraphics[width=0.45\columnwidth]{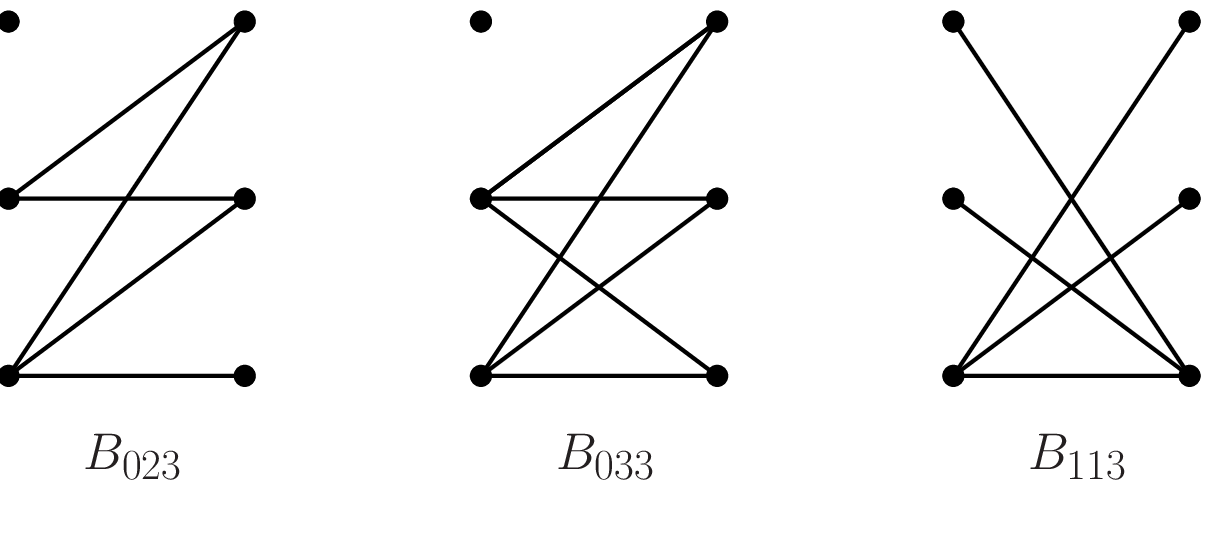}  
\caption{The graphs $B$ with $e(B) \ge 5$ and no perfect matching}
\label{fig:b033}
\end{center}
\end{figure}
We call the vertices of degree $3$ in $B_{113}$ the \emph{base vertices} of $B_{113}$ and the edge between them the \emph{base edge} of $B_{113}$.

The proof of the non-extremal case consists of four main steps.

\medskip

\noindent
{\bf Step 1:} We prove that for all but a constant number of vertices $x \in V(H) \backslash V(M)$, almost all
pairs $EF \in \binom{M}{2}$ are such that $L_x (EF) \cong B_{113}$. (See Claims~\ref{claim1}--\ref{claim4}.)

\noindent
{\bf Step 2:} We then show that this implies that $M$ must have size `close' to $d$ (see Claim~\ref{approxM}).

\noindent
{\bf Step 3:} Using Step~1, we show that there are 10 vertices $v_1,\dots , v_{10} \in V(H) \backslash V(M)$,
such that for almost all pairs $EF \in \binom{M}{2}$ not only does 
$L_{v_1} (EF)=\dots=L_{v_{10}} (EF) \cong B_{113}$ but further, for each such pair $EF$ the \emph{same} vertex $x$ plays the role
of the base vertex in $E$ (and the analogous statement holds for $F$ also). (See Claim~\ref{coro2} for the precise statement.)

\noindent
{\bf Step 4:} The information obtained in Steps~2 and~3 is then used to conclude that $H$ is `close' to $H_{n,d}$
(see Section~\ref{closesec}).

\medskip

To see how Fact~\ref{balanced} can be used in Step~1, suppose for example that $x_1$, $x_2$ and $x_3$ are unmatched vertices, that $E$ and $F$ 
are edges in $M$ and that the link graphs $L_{x_i}(EF)$ are identical (call this graph $B$).
The minimum degree condition implies that, for almost all unmatched vertices $x$, we have $e(L_{x}(EF)) \ge 5$. So let us assume this holds 
for $x_1, x_2, x_3$.
If $B$ contains a perfect matching, it is easy to see that we can transform $M$ into a (larger) matching which also 
covers the $x_i$, a contradiction.
If $B\cong B_{023},B_{033}$, we need to consider link graphs involving more than 2 edges from $M$
in order to obtain a contradiction. If $B=B_{113}$, we can use this to prove that we are `closer' to 
$H_{n,d}$. In particular, note that if $H=H_{n,d}$, then in the above example we have $B= B_{113}$.

To find a matching which is larger than $M$, we will often need several vertices whose link graphs 
with respect to some set of matching edges are identical (as in the above example).
We can usually achieve this with a simple application of the pigeonhole principle. 
But for this to work, we need to be able to assume that the number of vertices not covered by $M$ is fairly large.
This may not be true if e.g. we are seeking a perfect matching.
To overcome this problem, we apply the `absorbing method' which was first introduced in~\cite{rrs}.
The method (as used in~\cite{hps}) guarantees the existence of a small matching $M^*$ which can `absorb' any (very) small set of leftover vertices
$V'$ into a matching covering all of $V' \cup V(M^*)$. (The existence of $M^*$ is shown using a probabilistic argument.)
So if we are seeking e.g.~a perfect matching, it suffices to prove the existence of an almost perfect one outside $M^*$.
In particular, we can always assume that the set of vertices not covered by $M$ is reasonably large, as otherwise
we are done by the following lemma.
\begin{lemma}[H\`an, Person and Schacht~\cite{hps}]\label{absorb}

Given any $\gamma >0$ there exists an integer $n_0=n_0 (\gamma)$ such that the following holds. Suppose that $H$ is a $3$-uniform hypergraph
on $n\geq n_0$ vertices such that $\delta _1 (H) \geq (1/2+2\gamma)\binom{n}{2}$. 
Then there is a matching $M^*$ in $H$ of size $|M^*| \leq \gamma ^3n/3$ such 
that for every set $V' \subseteq V(H) \backslash V(M^*)$ with $\gamma ^6 n \geq |V'| \in 3 \mathbb Z$ there is a matching in $H$ covering 
precisely the vertices in $V(M^*) \cup V'$.
\end{lemma}

\section{Extremal case}\label{extremal}
The aim of this section is to show that hypergraphs which satisfy the degree condition in Theorem~\ref{generalmatch}
and are close to $H_{n,d}$ contain a $d$-matching.
\begin{lemma}\label{extremallemma}
There exist $\eps>0$ and $n_0 \in \mathbb N$ such that the following holds. Suppose that
$H$ is a $3$-uniform hypergraph on $n \geq n_0$ vertices and $d\le n/3$ is an integer. If
\begin{itemize}
\item $\delta _1 (H) > \binom{n-1}{2}- \binom{n-d}{2}$ and
\item $H$ is $\eps$-close to $H_{n,d}$,
\end{itemize}
then $H$ contains a $d$-matching.
\end{lemma}
We will first prove the lemma in the case when $H$ is not only close to $H_{n,d}$, but when for every vertex $v$
most of the edges of $H_{n,d}$ incident to $v$ also lie in $H$. More precisely,
given $\alpha >0$ and a $3$-uniform hypergraph $H$ on the same vertex set $V(H)$ as $H_{n,d}$,
we say that a vertex $v \in V(H)$ is 
\emph{$\alpha$-bad} if $|N_{H_{n,d}} (v) \backslash N_H (v)|> \alpha n^2$.
Otherwise we say that $v$ is \emph{$\alpha$-good}.
So if $v$ is $\alpha$-good then all but at most $\alpha n^2$ of the edges incident to $v$ in $H_{n,d}$
also lie in $H$.  We will now show that if $d\ge n/150$ then any such $H$ contains a $d$-matching.

\begin{lemma}\label{pmH} 
Let $0<\alpha<10^{-6}$ and let $n,d\in\mathbb{N}$ be such that $n/150\le d\le n/3$.
Suppose that $H$ is a $3$-uniform hypergraph on the same vertex set as $H_{n,d}$
and every vertex of $H$ is $\alpha$-good. Then $H$ contains a $d$-matching.
\end{lemma}
\proof Let $V$ and $W$ denote the vertex classes of $H_{n,d}$ of sizes $n-d$ and $d$
respectively. Consider the largest matching $M$ in $H$ which consists entirely of edges of type $ VVW$.
Let $V'$ denote the set of vertices in $V$ uncovered by~$M$. Define $W'$ similarly.
For a contradiction we assume that $|M|<d$. First note that $|M|\ge n/4$.
Indeed, to see this consider any vertex $w\in W'$. Since $w$ is $\alpha$-good but
$N_H(w)\cap \binom{V'}{2}=\emptyset$, it follows that%
     \COMMENT{otherwise have at least $\binom{2\sqrt{\alpha}n}{2}> \alpha n^2$
edges of $H_{n,d}$ which don't belong to $H$}
$|V'|\le 2\sqrt{\alpha}n$. Thus%
     \COMMENT{Note that if $d$ is not approx $n/3$ then this case cannot hold, ie $M$ will cover
all of $W$ and thus will be of size $d$.}
$|M|= |V\setminus V'|/2\ge (n-d-2\sqrt{\alpha}n)/2\ge n/4$.

Consider $v_1, v_2 \in V'$ and $w \in W'$ where $v_1 \not = v_2$. Given a pair $E_1 E_2$ of distinct
matching edges from $M$, we say that
$E_1 E_2$ is \emph{good for $v_1 v_2 w$} if there are all possible edges $E$ in $H$ which take the following form:
$E$ has type $VVW$ and contains one vertex from $\{v_1,v_2,w\}$, one vertex from $E_1$ and
one vertex from $E_2$. Note that if $E_1 E_2$ is good for $v_1 v_2 w$ then $H$ has a $3$-matching
which consists of edges of type $VVW$ and contains precisely the vertices in $E_1$, $E_2$ and $\{v_1,v_2,w\}$.
So if such a pair $E_1 E_2$ exists, we obtain a matching in $H$ that is larger than $M$, yielding a contradiction.

Since $|M|\geq n/4$ we have at least $\binom{n/4}{2} > n^2/40$ pairs of distinct matching edges
$E_1, E_2 \in M$.  Since $v_1, v_2$ and $w$ are $\alpha$-good there are at most $3 \alpha n^2 <n^2 /40$
such pairs $E_1 E_2$ that are not good for $v_1 v_2 w$.
So one such pair must be good for $v_1 v_2 w$, a contradiction.
\endproof

We now use Lemma~\ref{pmH} to prove Lemma~\ref{extremallemma}. Our strategy is to obtain a `small' matching
$M$ in $H$ that covers all `bad' vertices in $H$. We will construct $M$ in stages so as to ensure that
$H- V(M)$ satisfies the hypothesis of Lemma~\ref{pmH}. Thus we obtain
a $(d-|M|)$-matching $M'$ of $H- V(M)$, and hence a $d$-matching $M \cup M'$ of $H$.

\medskip

\noindent
{\bf Proof of Lemma~\ref{extremallemma}.}
Let $0<1/n_0\ll \eps \ll \eps ' \ll \eps''\ll \eps'''\ll 1$. By Theorem~\ref{bde} we may
assume that $d\ge n/100$. Suppose that $H$ is as in the
statement of the lemma and let $V$ and $W$ denote the vertex classes of $H$ of sizes $n-d$ and $d$ respectively.
Since $H$ is $\eps$-close to $H_{n,d}$, all but at most $3\sqrt{\eps}n$
vertices in $H$ are $\sqrt{\eps}$-good. Let $V^{bad}$ denote the set of $\sqrt{\eps}$-bad vertices in $V$. Define $W^{bad}$ similarly. 
So $|V^{bad}|,|W^{bad}|\leq 3\sqrt{\eps}n$. 

Define $c:= |W^{bad}|$, $V_1:=V\cup W^{bad}$ and $W_1:=W\backslash W^{bad}$. Thus $a:=|V_1|=n-d+c$ and
$b:=|W_1|=d-c$. Moreover,
$$\delta _1 (H[V_1])\geq \delta _1 (H)- \binom{b}{2}-(a-1)b>\binom{n-1}{2} - \binom{n-d}{2}- \binom{b}{2}-(a-1)b.$$
But $\binom{n-1}{2}= \binom{a-1}{2}+(a-1)b+\binom{b}{2}$ and so
$$\delta _1 (H[V_1]) > \binom{a-1}{2}-\binom{n-d}{2}=\binom{a-1}{2}-\binom{a-c}{2}.$$
Since  $c\leq 3\sqrt{\eps}n$ we can apply 
Theorem~\ref{bde} to obtain a matching $M_1$ of size $c$ in $H[V_1]$.

Let $H_1:=H- V(M_1)$ and $V_2:=V_1\backslash V(M_1)$. (Note that if $W^{bad} = \emptyset$ then $H_1=H$.) 
So $H_1$ has vertex classes $V_2$ and $W_1$ where $|V_2|=a-3c$.
Since $H$ is $\eps$-close to $H_{n,d} (V,W)$ and $3c\leq 9\sqrt{\eps}n  \ll \eps 'n$ we have that $H_1$ is
$\eps'$-close to $H_{|H_1|,b}(V_2, W_1)$.
By definition of $W_1$  all vertices in $W_1$ are $\eps'$-good in $H_1$.
Furthermore, if a vertex $v \in V(H_1)$ is $\eps'$-bad in $H_1$ then
$v \in V_2$ and $v \in V^{bad} \cup W^{bad}$. Let $V^{bad} _2$ denote the set of such vertices.
So $|V^{bad} _2|\leq 3 \sqrt{\eps}n$. If $V^{bad} _2= \emptyset$ then we can apply Lemma~\ref{pmH} to
obtain%
    \COMMENT{Need that $\frac{n-3c}{150}\le b=d-c\le \frac{n-3c}{3}$. But this is ok since
$n/100\le d\le n/3$.}
a $b$-matching $M_2$ in $H_1$. We thus obtain a matching
$M_1 \cup M_2$ of size $b+c=d$ in~$H$ . So we may assume that $V^{bad} _2\not = \emptyset$. 

We say that a vertex $v \in V^{bad}_2$ is \emph{useful} if there are at least
$\eps'n^2$ pairs of vertices $v'w \in V_2 W_1$ such that
$vv'w$ is an edge in $H_1$. Clearly we can greedily select a matching $M_2$ in $H_1$ such that
$m_2:=|M_2|\leq |V^{bad}_2|$ where $M_2$ covers
all useful vertices and consists entirely of edges of type $ V_2 V_2 W_1$.
Let $H_2:=H_1- V(M_2)$, $V_3:=V_2\backslash V(M_2)$ and $W_2:=W_1\backslash V(M_2)$.
Then $|V_3|=|V_2|-2m_2=a-3c-2m_2$ and $|W_2|=b-m_2$. Note that
\begin{align}\label{mindegcalc}
\delta_1(H)&>  \binom{n-1}{2} - \binom{n-d}{2} \ge
(1-\eps)\left(1-\left(1-\frac{d}{n}\right)^2\right)\frac{n^2}{2}\nonumber\\
& =(1-\eps)\left(\frac{2d}{n}-\frac{d^2}{n^2}\right)\frac{n^2}{2}
= (1-\eps)d\left(n-\frac{d}{2}\right).
\end{align}
Consider any vertex $v \in V^{bad} _2 \backslash V(M_2)$. Since $v$ is not useful, it must lie in more than
\begin{align*}
\delta_1(H)- & n|V(H)\setminus V(H_2)|-\eps'n^2-\binom{|W_2|}{2} \stackrel{(\ref{mindegcalc})}{\ge}
(1-\eps)d\left(n-\frac{d}{2}\right)-\eps' n^2-\eps' n^2-\frac{d^2}{2}\\
& \ge d(n-d)-\eps dn-2\eps' n^2\ge \frac{2dn}{3}-3\eps' n^2\ge 2\eps' n^2
\end{align*}
edges%
   \COMMENT{use that $|W_2|\le |W_1|\le d$}
of $H_2[V_3]$. Since $|V^{bad} _2|\leq 3 \sqrt{\eps}n$ we can greedily select a matching $M_3$ in $H_2[V_3]$ of
size $m_3:=|M_3|\leq |V^{bad} _2|$ which covers all the vertices in $H_2$ which lie in $V^{bad} _2$. 

Let $H_3:=H_2- V(M_3)$ and $V_4 := V_3 \backslash V(M_3)$. So $H_3$ has vertex classes 
$V_4$ and $W_2$ where $|V_4|=|V_3|-3m_3=a-3c-2m_2-3m_3$. Recall that every vertex in $V(H_1)\setminus V_2^{bad}$
is $\eps'$-good in $H_1$. Since $V_2^{bad}\subseteq V(M_2\cup M_3)$ and $|H_1|-|H_3|=3(|M_2|+|M_3|)\ll \eps' n$,
it follows that every vertex of $H_3$ is $\eps''$-good. So certainly for every vertex $w \in W_2$ there are at
least $|V_4||W_2|/2$ pairs $vw' \in V_4 W_2$ such that $vww'$ is an edge in $H_3$.
Thus we can greedily find a matching $M_4$ of size $m_3$ such that each edge in $M_4$ has type $V_4W_2W_2$. 

Let $H_4:=H_3- V(M_4)$, $V_5 := V_4 \backslash V(M_4)$ and $W_3 := W_2 \backslash V(M_4)$. So $H_4$ has  
vertex classes $V_5$ and $W_3$ of sizes $|V_5|=|V_4|-m_3=a-3c-2m_2-4m_3=n-d-2c-2m_2-4m_3$ and
$|W_3|=|W_2|-2m_3=b-m_2-2m_3=d-c-m_2-2m_3$. Moreover, every vertex of $H_4$ is $\eps'''$-good.
Thus we can apply Lemma~\ref{pmH} to $H_4$ to obtain%
     \COMMENT{Need that $\frac{|H_4|}{150}\le |W_3|\le \frac{|H_4|}{3}$. The first inequality holds since
$|W_3|\ge d-\eps' n\ge n/100-\eps' n$. To check the 2nd inequality note that $|H_4|=n-3c-3m_2-6m_3$ and
$3|W_3|=3d-3c-3m_2-6m_3$.}
a $|W_3|$-matching $M_5$ in $H_4$. But then $M_1 \cup M_2 \cup M_3 \cup M_4 \cup M_5$ is a matching of size
$c+m_2+m_3+m_3+|W_3|=d$ in $H$, as desired.
\endproof
We remark that the only point in the proof of Theorem~\ref{generalmatch} where we need the full strength of the 
minimum degree condition is when we apply Theorem~\ref{bde} to find the matching~$M_1$ in the proof of
Lemma~\ref{extremallemma}.


\section{Proof of Theorem~\ref{generalmatch}}\label{mainsec}
\subsection{Preliminaries}
We first define constants satisfying
\begin{align}\label{hier}
0<1/n_0  \ll 1/C  \ll \gamma'' \ll \gamma '\ll \gamma  \ll \eps ' \ll\eps \ll \eta ' \ll \eta
\ll \alpha' \ll \alpha \ll \rho' \ll\rho \ll \tau \ll 1.
\end{align}
Let $H$ be a $3$-uniform hypergraph on $n \geq n_0$ vertices such that
\begin{align}\label{mindeg}
\delta _1 (H) > \binom{n-1}{2}-\binom{n-d}{2} \geq (1-\gamma ')d(n-d/2),
\end{align}
where $d$ is an integer such that%
    \COMMENT{In particular, if $d=n/3$ then $n$ must be divisible by $3$ since $d\in \mathbb{N}$}
$1\leq d \leq n/3$. (Note that the second inequality in (\ref{mindeg}) follows from the same argument as~(\ref{mindegcalc}).)
 We wish to find a $d$-matching in $H$. Note that Theorem~\ref{bde} covers the
case when $d \leq n/100$. So we may assume that $n/100 \leq d \leq n/3$. 

Suppose $d \geq n/3 - \tau n$. Since $\tau \ll 1$, (\ref{mindeg}) gives us that
$ \delta _1 (H) \geq (1/2 + 2 \gamma '') \binom{n}{2}$.  
So by Lemma~\ref{absorb} there is a matching $M^*$ in $H$ of size $|M^*|\leq (\gamma'') ^3 n/3$ such that
for every set $V' \subseteq V(H) \backslash V(M^*)$ with $(\gamma'') ^6 n \geq |V'| \in 3 \mathbb Z$
there is a matching in $H$ covering precisely the vertices in $V(M^*) \cup V'$. 
If $n/100 \leq d < n/3 - \tau n$ we set $M^* := \emptyset$.

In both cases we define $H':= H- V(M^*)$. (So $H'=H$ if $n/100 \leq d < n/3 - \tau n$.)
Thus
\begin{align}\label{mindeg2}
\delta _1 (H') \ge \delta_1(H)-\gamma' n^2.
\end{align}
Let $M$ be the largest matching in $H'$. Clearly we may assume that $|M|<d$. Theorem~\ref{bde} implies that 
\begin{align}\label{Mbound}
  n/200 \leq |M| <d. 
\end{align}
Let $V_M:= V(M)$ and $V_0:= V(H')\backslash V_M$. So $|V_0| \leq n-|V_M|.$
If $n/100 \leq d < n/3 - \tau n$ then $|V_0|> n-3d >3 \tau n$. Suppose 
$d \geq n/3 - \tau n$. 
If $|V_0| \leq (\gamma'') ^6 n$, then by definition of $M^*$, there is a matching $M'$ in
$H$ containing all but at most two vertices from  $V(M^*) \cup V_0$.
But then $M \cup M'$ is a matching in $H$ of size $\lfloor n/3 \rfloor \geq d$, as desired. So in both cases we may assume that
\begin{align}\label{v0}
(\gamma'') ^6 n \leq |V_0| \leq  n-|V_M|.
\end{align}

\subsection{Finding structure in the link graphs}
In this section we show that `most' of our link graphs $L_{v} (EF)$ with $v \in V_0 $ and $EF \in \binom{M}{2}$
are copies of $B_{113}$ (recall that $B_{113}$ was defined after Fact~\ref{balanced}). 
\begin{claim}\label{claim1} There does not exist $v_1v_2v_3 \in \binom{V_0}{3}$ and $EF \in \binom{M}{2}$ such that
\begin{itemize}
\item $L_{v_1} (EF)=L_{v_2} (EF)=L_{v_3} (EF)$ and
\item $L_{v_1} (EF)$ contains a perfect matching.
\end{itemize}
\end{claim}
\proof The proof is identical to the proof of Fact~17 in~\cite{hps}. We include it here for completeness. 
Let $E=\{x_1, x_2, x_3\}$ and $F=\{y_1,y_2,y_3\}$ and suppose 
$x_1y_1$, $x_2y_2$ and $x_3y_3$ is a  perfect matching in $L_{v_1} (EF)$. Since these edges lie in 
$L_{v_i} (EF)$ for each $1 \leq i \leq 3$ the edges $v_1x_1y_1$, $v_2x_2y_2$ and $v_3x_3y_3$ lie in $H'$.
Replacing $E$ and $F$ in $M$ with these edges we obtain a larger matching in $H'$, a contradiction.
\endproof
We will now use Claim~\ref{claim1} to show that only a constant number of vertices $v \in V_0$
have `many' link graphs $L_{v} (EF)$ containing perfect matchings.
\begin{claim}\label{claim2} Let $V'_0$ denote the set of all those vertices $v \in V_0$ for which
there are at least $\eps n^2$ pairs $EF \in \binom{M}{2}$
such that $L_{v} (EF)$ contains a perfect matching. Then $|V'_0|\leq C$.
\end{claim}
\proof Let $G$ be the bipartite graph with vertex classes $V'_0$ and $\binom{M}{2}$
where $\{v, EF\}$ is an edge in $G$ precisely when 
$L_{v} (EF)$ contains a perfect matching. So $G$ contains at least $|V'_0|\eps n^2$ edges.
If $|V'_0|\geq C$ then there is a pair $EF \in \binom{M}{2}$
such that $d_G (EF) \geq C \eps \geq 3 \cdot 2^9$ (since $1/C \ll \eps$). Since there are $2^9$
labelled bipartite graphs with vertex classes $E$ and $F$,
there are $3$ vertices $v_1, v_2, v_3 \in V'_0$ such that  $L_{v_1} (EF)=L_{v_2} (EF)=L_{v_3} (EF)$ and
$L_{v_1} (EF)$ contains a perfect matching. This contradicts Claim~\ref{claim1}, as required.
\endproof
\begin{claim}\label{claim3}
Let $V''_0$ denote the set of all those vertices $v \in V_0$ for which there are at least $\eps n^2$ pairs
$EF \in \binom{M}{2}$ such that $L_{v} (EF)\cong B_{023}, B_{033}$. Then $|V''_0|\leq C$.
\end{claim}
\proof Suppose for a contradiction that $|V''_0|>C$. Given any $v \in V''_0$, define an auxiliary
oriented graph $G_v$ as follows: The vertex set of $G_v$ is $M$ and given $EF \in \binom{M}{2}$ there is
an edge directed from $E$ to $F$ precisely when $L_{v} (EF)\cong B_{023}, B_{033}$ where
$E$ is the vertex class that contains the isolated vertex in $L_{v} (EF)$. Since $v \in V''_0$, we have
that $e(G_v) \geq \eps n^2$.

We call a path $E_1\dots E_5$ of length $4$ in $G_v$ \emph{suitable} if its (directed) edges are 
$E_1 E_2, E_3 E_2 , E_3 E_4$ and  $E_5 E_4$. Our first aim is to find at least $\eps ' n^5$ suitable paths in $G_v$. 
Choose a partition $V_1, V_2$ of $V(G_v)$ such that 
$e _{G_v} (V_1, V_2) \geq e(G_v)/5 \geq \eps n^2/5$. (To see the existence of such a partition, consider the expected 
number of edges from $V_1$ to $V_2$ in a random
partition of $V(G_v)$.) Let $G'_v$ denote the undirected bipartite graph with vertex classes $V_1$
and $V_2$ whose edges are all those edges in $G_v$ that are oriented from $V_1$ to $V_2$.
Since $e(G'_v) \geq \eps n^2/5$, $G'_v$ contains a subgraph $G''_v$ with $\delta (G''_v) \geq d(G'_v)/2 \geq \eps n/5$.
Thus we can greedily find at least 
$$\frac{1}{2} \cdot \frac{\eps n}{5}\left( \frac{\eps n}{5}-1\right) \dots
\left( \frac{\eps n}{5}-4\right) \geq \eps ' n^5$$
paths of length $4$ in $G''_v$ whose endpoints both lie in $V_1$. By definition of $G''_v$, each of
these paths corresponds to a suitable path in $G_v$.

Consider a suitable path $E_1\dots E_5$ in $G_v$. So
$L_v (E_2 E_3), L_v (E_3 E_4) \cong B_{023}, B_{033}$ with the isolated vertex in both graphs lying in $E_3$.
Choose edges $e_1$ of $L_v(E_2E_3)$ and $e_2$ of $L_v(E_3 E_4)$ such that $e_1$ and $e_2$ are disjoint. Since 
$L_v (E_1 E_2) \cong B_{023}, B_{033}$ and $E_1$ contains the isolated vertex in this graph,
there is a $2$-matching $\{e_3, e_4\}$ in 
$L_v (E_1 E_2)$ that is disjoint from $e_1$. Similarly since $L_v (E_4 E_5) \cong B_{023}, B_{033}$ 
and $E_5$ contains the isolated vertex in this graph, there is a $2$-matching $\{e_5, e_6\}$ in 
$L_v (E_4 E_5)$ that is disjoint from $ e_2$. Hence $L_v (E_1 E_2E_3E_4E_5)$ contains a $6$-matching
$\{e_1, e_2, e_3, e_4,e_5, e_6\}$.

Let $G$ be the bipartite graph with vertex classes $V''_0$ and the set $(M)^5$ of all ordered $5$-tuples of elements of $M$ where
$\{v, E_1E_2E_3E_4E_5\}$ is an edge in $G$ precisely when $E_1\dots E_5$ is a suitable path  in $G_v$.
So $G$ contains at least $|V''_0|\eps ' n^5 $ edges. 

Since $|V''_0| >C$ there exists
$E_1E_2E_3E_4E_5 \in (M)^5$  such that $d_G (E_1E_2E_3E_4E_5) \geq C\eps ' \geq 6 \cdot 2^{36}$. Further, 
there are at most $2^{36}$ distinct
graphs in the collection of all those graphs $L_v(E_1E_2E_3E_4E_5)$ for which $v \in N_G (E_1E_2E_3E_4E_5)$.
Thus there are 6 vertices $v_1,\dots , v_6 \in V''_0$ such that $v_1, \dots ,v_6 \in N_G (E_1E_2E_3E_4E_5)$ and 
$L_{v_1} (E_1E_2E_3E_4E_5)= \dots = L_{v_6} (E_1E_2E_3E_4E_5)$. Let $\{x_1 y_1, \dots , x_6 y_6\}$ be a $6$-matching in 
$L_{v_1} (E_1E_2E_3E_4E_5)$. So  $\{v_1 x_1 y_1, \dots , v_6 x_6 y_6\}$ is a $6$-matching in $H'$. Replacing the edges
$E_1,\dots, E_5$ in $M$ with $\{v_1 x_1 y_1, \dots , v_6 x_6 y_6\}$ we obtain a larger matching, a contradiction.
\endproof

\begin{claim} \label{claima}
Let $V'''_0$ denote the set of all those vertices $v \in V_0$ which fail to satisfy
\begin{equation} \label{bipbound}
 e(L_v(V_0,V_M)) \le ( 1+\sqrt{\gamma'}) |V_0||M|.
\end{equation}
Then $|V'''_0|\leq C$.
\end{claim}
\proof
Suppose for a contradiction that $|V'''_0|> C \ge 2/\gamma'$.
Given an edge $E$ in $M$, we say that $E$ is \emph{good for $v \in V'''_0$} 
if at least two vertices in $E$ have degree at least $3$ in $L_v(E,V_0)$. 
For every $v \in V'''_0$, there are at least $\gamma' |M|$  edges in $M$ 
which are good for $v$.
(To see this, suppose there are fewer edges which are good for $v$.
Then 
\begin{align*}
e(L_v(V_0,V_M)) & < (1-\gamma')|M|(4+|V_0|)+\gamma'|M| \cdot 3|V_0|\\
& \le |M||V_0|\left( (1-\gamma')(1+\gamma')+3\gamma'\right)
 \le (1+\sqrt{\gamma'}) |V_0||M|,
\end{align*}
a contradiction to the fact that $v \in V'''_0$.)
This in turn implies that there are  $v_1,v_2 \in V'''_0$ and an edge $E$ in $M$ which is good 
for both $v_1$ and $v_2$.
Then the definition of `good' implies that are disjoint edges $e_1 \in L_{v_1}(E,V_0)$ and $e_2 \in L_{v_2}(E,V_0)$ 
which do not contain $v_1$ or $v_2$.
Now we can enlarge $M$ by removing $E$ and adding $v_1e_1$ and $v_2e_2$.
This contradiction to the maximality of $M$ proves the claim. \endproof

\begin{claim}\label{claimb}
Every vertex $v \in V _0 \backslash V'''_0$ satisfies
\begin{equation*}
e(L_v(V_M)) \ge \left(5-\gamma \right) \binom{|M|}{2} .
\end{equation*}
\end{claim}
\proof
Suppose $v \in V_0 \backslash V'''_0$. 
Then as $e(L_v(V_0))=0$
\begin{eqnarray*}
e(L_v(V_M)) & \stackrel{(\ref{mindeg2})}{\ge} & \delta_1(H) - e(L_v(V_0,V_M)) -\gamma ' n^2\\
& \stackrel{(\ref{mindeg}),(\ref{bipbound})}{\ge} & (1-\gamma ')d(n-d/2) - \left(1+\sqrt{\gamma'} \right) |V_0||M|-\gamma ' n^2.
\end{eqnarray*}
Now note that the function $d(n-d/2)$ is increasing in $d$ for $d \le n/3$. So
\begin{align*}
e(L_v(V_M)) & \ge (1-\gamma ')|M|\left( n-\frac{|M|}{2} \right) - \left(1+\sqrt{\gamma'} \right) (n-3|M|)|M|-\gamma ' n^2 \\
& \ge \left( n|M|-\frac{|M|^2}{2} -\gamma' n|M|\right)- \left( n|M|-3|M|^2+\sqrt{\gamma'}n|M| \right)-\gamma ' n^2 \\
& \stackrel{(\ref{Mbound})}{\ge} \frac{5|M|^2}{2}-400\sqrt{\gamma'} |M|^2\ge (5-\gamma) \binom{|M|}{2},
\end{align*}
which completes the proof of the claim.
\endproof

\begin{claim}\label{claim4}
Let $V''''_0$ denote the set of all those vertices $v \in V_0 \backslash V'''_0$ for which there are at least $\eta n^2$ pairs
$EF \in \binom{M}{2}$ such that $L_{v} (EF)$ contains at most $4$ edges. Then $|V''''_0|\leq 2C$.
\end{claim}
\proof Suppose for a contradiction that $|V''''_0| >2C$. Let $v \in V''''_0$. 
At most $3|M|$ edges $vv_1v_2$ in $H$ containing $v$ are such that
$v_1$ and $v_2$ lie in the same edge $E \in M$. 
Thus Claim~\ref{claimb} implies that
\begin{align}\label{sumcon}
\sum _{EF \in\binom{M}{2}} e(L_v (EF))\geq (5-\gamma) \binom{|M|}{2}-3|M|\geq 5\binom{|M|}{2}- \gamma n^2.
\end{align}
Let $c$ denote the number of pairs $EF \in \binom{M}{2}$ such that $L_v (EF)$ contains at most $4$ edges.
Then $ c \geq \eta n^2$ and so (\ref{sumcon}) implies that there are at least $\eta' n^2$ pairs
$EF \in \binom{M}{2}$ such that $L_v (EF)$ contains at least $6$ edges. Indeed, suppose
that this is not the case. Then
\begin{align*}
\sum _{EF \in \binom{M}{2}} e(L_v (EF))  & \leq 4c + 9\eta ' n^2 + 5\left[\binom{|M|}{2}-c\right] =
5 \binom{|M|}{2} -c + 9 \eta' n^2\\ & < 5 \binom{|M|}{2}- \gamma n^2
\end{align*}
since $\gamma  \ll \eta ' \ll \eta $. This contradicts (\ref{sumcon}), as desired.

Recall from Fact~\ref{balanced} 
that a balanced bipartite graph $B$ on $6$ vertices that contains at least $6$ edges either
has a perfect matching or $B \cong B_{033}$. Thus, given any $v \in V''''_0$
there are at least $r \geq  \eta ' n^2 /2 \geq \eps n^2$ pairs%
\COMMENT{this is why we need $\eta \gg \eps$}
$E_1F_1, \dots , E_r F_r \in \binom{M}{2}$
such that  either
\begin{itemize} 
\item $L_v (E_i F_i)$ contains a perfect matching for all $1\leq i\leq r$ or,
\item $L_v (E_i F_i)\cong  B_{033}$ for all $1\leq i\leq r$.
\end{itemize}
So since $|V''''_0|>2C$ one of the following holds:
\begin{itemize} 
\item[($\alpha_1$)] There are more than $C$ vertices $v \in V''''_0$ for which there are at
least $\eps n^2$ pairs $EF\in \binom{M}{2}$ such that $L_v (EF)$ contains a perfect matching.
\item[($\alpha _2$)] There are more than $C$ vertices $v \in V''''_0$ for which there
are at least $\eps n^2$ pairs $EF\in \binom{M}{2}$ such that $L_v (EF)\cong  B_{033}$.
\end{itemize}
In either case we get a contradiction: ($\alpha _1$) contradicts Claim~\ref{claim2} and
($\alpha _2$) contradicts Claim~\ref{claim3}.
\endproof
Recall from Fact~\ref{balanced} that if $B$ is a balanced bipartite graph on $6$ vertices with $e(B)=5$ then either $B$ contains
a perfect matching or $B \cong B_{023},B_{113}$. If $e(B)\geq 6$ then either $B$ contains a perfect
matching or $B \cong B_{033}$. Thus Claims~\ref{claim2},~\ref{claim3},~\ref{claima} and~\ref{claim4} together
imply that all vertices $v \in V_0 \setminus (V_0' \cup V_0'' \cup V_0''' \cup V_0'''')$ satisfy%
\COMMENT{uglier but formally necessary as we use other properties and not just $(\beta$)}
\begin{itemize}
\item[($\beta$)] $L_v (EF) \cong B_{113}$ for at least $\binom{|M|}{2}-2\eps n^2 -\eta n^2
\geq (1-\alpha') \binom{|M|}{2}$ pairs $EF \in \binom{M}{2}$.
\end{itemize}
Let $V^* _0:=V_0\setminus (V_0' \cup V_0'' \cup V_0''' \cup V_0'''')$.
Thus 
$$
|V_0 \setminus V_0^*| \le 5C.
$$ 
Moreover, each $v\in V_0^*$ satisfies
\begin{align}\label{nicebound}
e(L_v(V_M))\le 5(1- \alpha ')\binom{|M|}{2}+9 \alpha ' \binom{|M|}{2}+3|M|\le 5(1+\alpha')\binom{|M|}{2}.
\end{align}
Here the term $3|M|$ accounts for the edges which have both endpoints in the same matching edge of $M$.

We can now show that $M$ has almost the required size. (This corresponds to Step~2 in the proof outline.)
This will be used in Section~\ref{closesec} to prove that $H$ is close to $H_{n,d}$.
\begin{claim}$|M|>d- \alpha n$. \label{approxM}
\end{claim}
\proof Assume for a contradiction that $|M|\leq d- \alpha n$. Consider any $v \in V^* _0$.
Then
\begin{equation}\label{eqdegv}
d_{H'} (v) \stackrel{(\ref{mindeg}), (\ref{mindeg2})}{\geq} (1- \gamma ')d(n-d/2)-\gamma ' n^2\ge d(n-d/2)- 2 \gamma ' n^2.
\end{equation}
Also $e(L_v(V_0))=0$ since $M$ is maximal. Thus
\begin{eqnarray*}
d_{H'} (v) &= & e(L_v (V_M))+e(L_v (V_0,V_M))
\stackrel{(\ref{bipbound}),(\ref{nicebound})}{\leq} 5(1+\alpha')\binom{|M|}{2}+(1+ \sqrt{\gamma'})|V_0||M|\\
& \leq & 5(1+\alpha')\binom{|M|}{2} +\left( |M|(n-3|M|)+ \sqrt{\gamma '} n^2 \right)\\
& \le & |M|(n-|M|/2)+\sqrt{\alpha '}n^2
< (d-\alpha n)(n-d/2+\alpha n/2)+\sqrt{\alpha '}n^2 \\
& < & d(n-d/2) -2\gamma' n^2,
\end{eqnarray*}
a contradiction to~(\ref{eqdegv}), as desired. (In the third line we again used that the function $d(n-d/2)$ is increasing in $d$
for $d\le n/3$.)
\endproof
In the next sequence of claims, we will show that there are vertices $v_1,\dots,v_{10} \in V_0^*$ whose 
 link graphs $L_{v_i}(V_M)$ are very similar to each other
(see Claim~\ref{coro2} for the precise statement). (This corresponds to Step~3 in the proof outline.)
\begin{claim}\label{claim5} Suppose $v_1, \dots, v_{10}\in V^* _0$ are distinct vertices such that for some $EF \in \binom{M}{2}$,
$L_{v_1} (EF), \dots, L_{v_{10}} (EF) \cong B_{113}$. Then $L_{v_1} (EF)= \dots =L_{v_{10}}(EF)$.
\end{claim}
\proof We suppose for a contradiction that the claim does not hold. Since there are $9$ labelled bipartite graphs
with vertex classes $E$ and $F$ which
are isomorphic to $B_{113}$, two of the $L_{v_i} (EF)$ must be the same. So we may assume that
$L_{v_1} (EF)= L_{v_2} (EF)$ but $L_{v_1} (EF)\not = L_{v_3} (EF)$.
Let $E=\{x_1, x_2, x_3\}$ and $F=\{y_1,y_2,y_3\}$.
Suppose $E(L_{v_1} (EF))=E(L_{v_2} (EF))=\{x_1 y_1, x_1y_2, x_1 y_3, x_2 y_1, x_3 y_1\}$. (So $x_1 y_1$ is the base
edge of $L_{v_1} (EF)$ and $L_{v_2} (EF)$ as defined after Fact~\ref{balanced}.)
Since $L_{v_1} (EF)\not = L_{v_3} (EF)$ there is an edge $e \in L_{v_3} (EF) \backslash L_{v_1} (EF)$.
We may assume $e=x_3 y_3$. Replacing $E$ and $F$ with
$v_1 x_1 y_2, v_2 x_2 y_1$ and $ v_3 x_3 y_3$ in $M$ we obtain a larger matching, a contradiction.
\endproof
Choose distinct $v_1, \dots , v_{10}\in V^* _0$ which will be fixed throughout the remainder of the proof.
\begin{claim}\label{claim6}
There is a set $\E$ of at least $(1- \alpha)|M|$ matching edges $E \in M$ such that for each $E\in \E$
there are at least $(1-\alpha)|M|$ edges $F \in M$
for which
$$L_{v_1} (EF)= \dots = L_{v_{10}} (EF) \cong B_{113}.$$
\end{claim}
\proof
By ($\beta$) and Claim~\ref{claim5} there are at least $(1-10\alpha') \binom{|M|}{2}$
pairs $EF \in \binom{M}{2}$ such that $L_{v_1} (EF)= \dots = L_{v_{10}} (EF) \cong B_{113}$.
This in turn immediately implies the claim.
\endproof

\begin{claim}\label{claim7}
For every $E \in \E$ there is a set $\F_E$ of at least $(1- 2 \alpha)|M|$ edges in $M$ such that
\begin{itemize}
\item[($\delta _1$)] $L_{v_1} (EF)= \dots = L_{v_{10}} (EF) \cong B_{113}$ for each $F\in \F_E$ and
\item[($\delta _2$)] in each of the $L_{v_1} (EF)$ with $F\in \F_E$ the same  vertex $x$ plays
the role of the base vertex in $E$. (Recall that the base vertices of $B_{113}$ are the vertices of degree $3$.)
\end{itemize}
\end{claim}
\proof Since $E\in \E$ there is a set $\F'_E$ of at least $(1- \alpha)|M|$ edges in $M$ such that
$L_{v_1} (EF)= \dots = L_{v_{10}} (EF) \cong B_{113}$ for each $F\in \F'_E$. 
Let $\F_E:=\F'_E\cap \E$. Then $|\F_E|\ge (1- 2 \alpha)|M|$ and 
for each $F\in \F_E$ there are at least $(1-\alpha)|M|$ edges $F' \in M$ for which
$L_{v_1} (FF')= \dots = L_{v_{5}} (FF') \cong B_{113}$.

We claim that $\F_E$ satisfies the claim. Certainly $\F_E$ satisfies ($\delta _1$).
Suppose for a contradiction that there are $F_1,F_2\in \F_E$ such that the vertex $x_1 \in E$
that plays the role of a base vertex in $L_{v_1} (EF_1)$ is different from the vertex
$x_2 \in E$ that plays the role of a base vertex in $L_{v_1} (EF_2)$. Let $F' \in M$ be such that 
$L_{v_1} (F_2F')= \dots = L_{v_5} (F_2F') \cong B_{113}$, and $F'\not = E,F_1$.

Since $L_{v_1} (EF_1)\cong B_{113}$ and $x_1 \not =x_2$, there exists a $2$-matching
$\{e_1, e_2\}$ in $L_{v_1} (EF_1)$ that is disjoint from $x_2$. 
Similarly since $L_{v_1} (F_2F')\cong B_{113}$ there exists a $2$-matching $\{e_3, e_4\}$ in $L_{v_1} (F_2F')$.
Since $x_2 \in E$ is a base vertex in $L_{v_1} (EF_2)$, there is an edge $e_5$ from $x_2$ to the vertex
in $F_2$ that is uncovered by $\{e_3, e_4\}$. So $\{e_1, e_2, e_3, e_4, e_5\}$ is a $5$-matching in
$L_{v_1} (F_1 E F_2 F')$.
We have chosen $F_1, F_2 $ and $F'$ so that 
$L_{v_1} (F_1 E F_2 F')=L_{v_2} (F_1 E F_2 F')= \dots = L_{v_5} (F_1 E F_2 F')$. Thus
$M':=\{v_1e_1, v_2e_2, v_3e_3, v_4 e_4, v_5 e_5\}$
is a $5$-matching in $H'$ that contains only vertices from
$E\cup F' \cup F_1 \cup F_2 \cup \{v_1, v_2, v_3, v_4, v_5\}$. Replacing $E, F', F_1$ and $F_2$ in $M$
with the edges in $M'$ yields a larger matching, a contradiction.
\endproof
 
Given $E \in \E$, we call the unique vertex $x \in V(E)$ satisfying $(\delta _2)$ a \emph{bottom vertex}. 
If $y \in E$ is such that $y \not = x$ then we say that $y$ is a \emph{top vertex}.
So each $E\in \E$ contains one bottom vertex and two top vertices whereas none of the at most
$\alpha |M|$ edges in $M\setminus \E$ contains a top or bottom vertex.
\begin{claim}\label{coro2}
There are at least $(1- 6{\alpha} )|M|^2 /2$ pairs $EF \in \binom{M}{2}$ such that
\begin{itemize}
\item[($\eps _1$)] $L_{v_1} (EF)= \dots = L_{v_{10}} (EF) \cong B_{113}$;
\item[($\eps _2$)] both $E$ and $F$ contain a bottom vertex $w$ and $z$ respectively;
\item[($\eps _3$)] $wz$ is the base edge of $L_{v_1} (EF)$.
\end{itemize}
\end{claim}
\proof 
Consider the directed graph $G$ whose vertex set is $M$ and in which there is a directed edge from
$E$ to $F$ if $E\in \E$ and $F\in \F_E$. Claims~\ref{claim6} and~\ref{claim7} together imply that
$G$ has at least $(1-3\alpha)|M|^2$ edges and thus at least $(1- 6{\alpha} )|M|^2 /2$
pairs $EF$ of vertices in $G$ must be joined by a double edge.%
   \COMMENT{since $2(1-6\alpha)|M|^2/2+ 6\alpha |M|^2/2\le (1-3\alpha)|M|^2$}
But each such pair $EF$ satisfies the claim.
\endproof

\subsection{Showing that $H$ is $\sqrt{\rho}$-close to $H_{n,d}$} \label{closesec}
We have now collected all the information we need for showing that $H$ is close to $H_{n,d}(V,W)$,
where $W$ will be constructed from the set of bottom vertices in $M$. More precisely,
let $W'$ denote the set of all the bottom vertices. So Claims~\ref{approxM} and~\ref{claim6} together imply that 
\begin{equation}\label{sizeW'}
d-2\alpha n \le (1-\alpha )|M| \leq |\E|= |W'| \leq |M| \leq d.
\end{equation}
Let $V'$ denote the set of all the top vertices in $H$. Thus 
\begin{equation}\label{sizeV'}
2d-4\alpha n \le 2(1- \alpha)|M|\leq |V'|=2|W'|\leq 2d.
\end{equation}
Choose a partition $V,W$ of $V(H)$ such that $|W|=d$, $W'\subseteq W$, $V'\subseteq V$.
Note that since (\ref{sizeW'}) implies that
$|W\setminus W'|\le 2 \alpha n$, all but at most $2 \alpha n$ vertices of $V_0$ lie in $V$.%
\COMMENT{was $3\alpha n$ before}
Our aim  is to show that $H$ is $\sqrt{\rho}$-close to $H_{n,d}(V,W)$.
Note that showing this proves Theorem~\ref{generalmatch} as we can apply Lemma~\ref{extremallemma} since we
chose $\rho \ll 1$ in (\ref{hier}).

\begin{claim}\label{claimV'V0V0}
$H$ does not contain an edge of type $V'V_0V_0$.
\end{claim}
\proof
Suppose that the claim is false and let $v'vv_0$ be an edge of $H$ with  $v'\in V'$ and $v,v_0\in V_0$.
Let $E \in \mathcal{E}$ be the matching edge containing $v'$.
Take any $F \in \mathcal{F}_E$. 
Take any 2 vertices from $v_1,\dots,v_{10}$ which are not equal to $v_0$ or $v$, call them $x$ and $y$.
Since $v'$ is a top vertex of $E$, it follows that $L_x(EF)$ contains a 2-matching
$e_1,e_2$ avoiding $v'$. Note that this is also a 2-matching in $L_y(EF)$.
Now we can enlarge $M$ by removing $E,F$ and adding $v'vv_0$, $xe_1$ and $ye_2$.
This contradicts the maximality of $M$ and proves the claim.
\endproof

\begin{claim}\label{claimlink} $\ $
\begin{itemize}
\item $H$ contains at least $(1-\rho')|W'||V'||V_0|$ edges of type $W'V'V_0$.
\item $H$ contains at least $(1-\rho')|V_0|\binom{|W'|}{2}$ edges of type $W'W'V_0$.
\item $H$ contains at most $\rho' |V_0|\binom{|V'|}{2}$ edges of type $V'V'V_0$.
\end{itemize}
\end{claim}
\proof
To see the first part of the claim, consider any $v\in V^*_0$ and any pair $w',v'$ with $w'\in W'$ and $v'\in V'$.
Both $w',v'$ could lie in the same matching edge from $M$, but there are at most $3|M|$ such pairs.
Also, $w',v'$ could lie in a pair $E,F$ of matching edges from $M$ for which either $L_v(EF) \not\cong B_{113}$
or which does not satisfy ($\eps_1$)--($\eps_3$) in Claim~\ref{coro2}.
But ($\beta$) and Claim~\ref{coro2} together imply that there are at most $\sqrt{\alpha} n^2$ such pairs $E,F$.
So suppose next that $w',v'$ lie in a pair $E,F$ satisfying $L_v(EF) \cong B_{113}$
and ($\eps_1$)--($\eps_3$). Then $L_v (EF), L_{v_1} (EF), \dots, L_{v_{9}} (EF) \cong B_{113}$
and so $L_{v} (EF)= L_{v_1} (EF)= \dots = L_{v_{9}} (EF)$ by Claim~\ref{claim5}. Conditions~$(\eps_2)$ and $(\eps_3)$
now imply
that $w'v'\in E(L_v(W',V'))$. So
$$
e(L_v(V',W')) \ge |V'||W'| - 2\sqrt{\alpha} n^2
\ge (1-\rho'/2) |V'||W'|.
$$
Summing over all vertices $v\in V^*_0$ and using that $|V_0\setminus V^*_0|\le 5C$
implies the first part of the claim. The remaining parts of the claim can be proved similarly.
\endproof

\begin{claim}\label{claimedges}
$H$ contains at least $|W'|\binom{|V_0|}{2}- \rho  n^3$ edges of type $W'V_0V_0$.
\end{claim}
\proof
Consider any $v \in V_0$.
By Claim~\ref{claimV'V0V0} there are no edges in $L_v(V(H))$ with one endpoint in $V'$ and the other in $V_0$.
By~(\ref{sizeW'}) there are at most $3\alpha |M|n \leq 3\alpha n^2$ edges in $L_v (V(H))$ with one endpoint
in $V_M\backslash (V' \cup W')$ and the other in $V_0$. Furthermore, $L_v (V_0)$ contains no edges.
Thus,
\begin{eqnarray*}
e(L_v(W',V_0)) & \ge & \delta_1(H') - e(L_v(V_M)) - 3{\alpha}n^2\\
& \stackrel{(\ref{mindeg}), (\ref{mindeg2}), (\ref{nicebound})}{\ge} & 
(1-\gamma ')d\left(n-\frac{d}{2}\right) -\gamma' n^2-5(1+\alpha ')\binom{|M|}{2}- 3{\alpha}n^2\\
& \stackrel{(\ref{Mbound})}{\ge} & (1-\gamma')|M|\left(n-\frac{|M|}{2}\right) -(5+\sqrt{\alpha})\frac{|M|^2}{2} \\
& \ge & |M|(n -3|M|)- \sqrt{\alpha}|M|n \ge|W'||V_0|- \rho'  n^2.
\end{eqnarray*} 
As earlier, here we use the fact that the function $d(n-d/2)$ is increasing in $d$ for $d\leq n/3$.
Summing over all vertices $v\in V^*_0$ and using the fact that $|V_0\setminus V^*_0|\le 5C$ now proves the claim. 
\endproof

\begin{claim}\label{claimtypes} $\ \ $
\begin{itemize}
\item $H$ contains at least $(1-\rho)|W'|\binom{|V'|}{2}$ edges of type $W'V'V'$.
\item $H$ contains at least $(1-\rho)|V'|\binom{|W'|}{2}$ edges of type $W'W'V'$.
\end{itemize}
\end{claim}
\proof
First note that the last part of Claim~\ref{claimlink} implies that all but at most $2\sqrt{\rho'} n$ vertices%
\COMMENT{each edge is counted at most twice}
$x\in V'$
lie in at most $\sqrt{\rho'} |V'||V_0|$ edges of type $V'V'V_0$. Call such vertices $x$ \emph{useful}.
Consider any useful $x$. Then $x \in E'$ for some $E' \in \E\subseteq M$. Further, since $x$ is a top vertex in $E'$,
certainly there exists an edge $F' \in M$
such that $L_{v_1} (E'F')= L_{v_2} (E'F')\cong B_{113}$, where $x$ is not a base vertex in 
$L_{v_1} (E'F')$. So $L_{v_1} (E'F') $ contains a $2$-matching $\{ e_1, e_2\}$ which avoids $x$. 

Consider any pair 
$EF \in \binom{M\backslash \{E',F'\}}{2}$ satisfying ($\eps _1$)--($\eps _3$). We claim that
$L_{x} (EF) \subseteq L_{v_1} (EF)$.
Indeed, if not then there exist disjoint edges $e_3, e_4$ and $e_5$ such that $e_3 \in E(L_x (EF))$ and
$e_4, e_5 \in E(L_{v_1} (EF))$.
Since $L_{v_1} (E'F')= L_{v_2} (E'F')$ and since $EF$ satisfies ($\eps _1$) we have that
$v_1e_1, v_2 e_2 , x e_3, v_3 e_4$ and $v_4 e_5$ are edges in $H'$.
Replacing $E,F,E',F'$ with $v_1e_1, v_2 e_2 , x e_3, v_3 e_4$ and $v_4 e_5$ in $M$ yields a larger matching
in $H'$, a contradiction. So indeed $L_{x} (EF) \subseteq L_{v_1} (EF)$.

There are at least $(1- 6{\alpha})|M|^2/2-2|M| \geq (1-7{\alpha})|M|^2/2$ pairs
$EF \in \binom{M\backslash \{E',F'\}}{2}$ 
satisfying ($\eps _1$)--($\eps _3$). We claim that at most $\rho^2 |M| ^2 /2$ of these pairs $EF$ are such that 
$L_x (EF)$ contains fewer than $5$ edges. Indeed, suppose not. Since for such $EF$,
$L_{x} (EF) \subseteq L_{v_1} (EF)\cong B_{113}$, the number of edges of $H$ which contain $x$ and
have no endpoint outside $V_M$ is at most
$$
4\cdot \rho  ^2 |M|^2 /2 + 5 \cdot (1-7\alpha  - \rho^2 ) |M|^2 /2+9 \cdot 7{\alpha } |M|^2/2 +3|M|
\le (5+30\alpha-\rho^2)|M|^2/2.
$$
Here the third term accounts for edges between pairs not satisfying ($\eps _1$)--($\eps _3$)
and the final term for edges with 2 vertices in the same matching edge from $M$. Let us now bound the number of
edges containing $x$ which have an endpoint outside $V_M$. There are at most $|W'|(n-3|M|)\leq|M|(n-3|M|)$ such edges
having an endpoint in $W'$ and at most $\sqrt{\alpha}n^2$ such edges having an endpoint outside $V'\cup W'\cup V_0$.
Since $H$ has no edge of type $V'V_0V_0$ by Claim~\ref{claimV'V0V0}, the only other such edges consist of $x$,
one vertex in $V'$ and one vertex in $V_0$. But since $x$ is useful the number of such edges is at most $\sqrt{\rho'} |V'||V_0|$.
Thus in total there are at most $|M|(n-3|M|)+2\sqrt{\rho'} n^2$ edges which contain $x$ and
have an endpoint outside $V_M$. So the degree of 
$x$ in $H$ is at most 
\begin{eqnarray*}
(5+30\alpha-\rho^2)|M|^2/2+|M|(n-3|M|)+ 2\sqrt{\rho'} n^2 & \le & |M|(n-|M|/2)-\rho^3 n^2\\
& {\le} & d(n-d/2)-\rho^3 n^2 \stackrel{(\ref{Mbound}),(\ref{mindeg})}{<} \delta _1 (H),
\end{eqnarray*}
a contradiction.
Thus there are at least $(1- 7{\alpha } - \rho^2 )|M|^2 /2$ pairs $EF \in \binom{M\backslash \{E',F'\}}{2}$
satisfying ($\eps _1$)--($\eps _3$) such that $L_{x} (EF) = L_{v_1} (EF)\cong B_{113}$.
Let $\cP$ denote the set of such pairs.

Now consider any pair $w',v'$ with $w'\in W'$ and $v'\in V'\setminus \{x\}$.
Both $w',v'$ could lie in the same matching edge from $M$, but there are at most $3|M|$ such pairs.
Also, $w',v'$ could lie in a pair $E,F$ of matching edges which does not belong to $\cP$.
But there at most $5\rho^2 |M|^2$ such pairs $w',v'$.
So suppose next that $w',v'$ lies in a pair $E,F$ belonging to $\cP$.
Since $L_{x} (EF)= L_{v_1} (EF) \cong B_{113}$ and $EF$ satisfies ($\eps _2$) and ($\eps _3$)
it follows that $w'v'\in E(L_{x} (EF))$. Thus $e(L_x(W',V'))\ge (1-6\rho^2)|W'||V'|$.
Summing over all useful vertices $x\in V'$ proves the first part of the claim. The second part follows
similarly (the only change is that we consider a pair $w_1',w_2' \in W'$ in the final paragraph).
\endproof

Claims~\ref{claimlink}--\ref{claimtypes} together with~(\ref{sizeW'}) and~(\ref{sizeV'}) now show that $H$ contains all
but at most $\sqrt{\rho} n^3$ edges of type $WVV$ and $WWV$ and thus $H$ is $\sqrt{\rho}$-close to
$H_{n,d}(V,W)$. Hence $H$ contains a perfect matching by Lemma~\ref{extremallemma}.
\\
\newline \noindent
{\bf Remark.} One can also obtain Theorem~\ref{generalmatch} by proving the result only in the case when $d=\lfloor n/3 \rfloor$. Indeed, 
suppose that $H$ is as in the theorem. Let $a:=\lfloor (n-3d)/2\rfloor$. Obtain a new $3$-uniform hypergraph $H'$
from $H$ by adding $a$ new vertices to $H$ such that each of these vertices forms an edge with all pairs of vertices in $H'$. It is not
hard to check that $\delta _1 (H') > \binom{|H'|-1}{2}-\binom{|H'|- \lfloor |H'|/3\rfloor}{2}$ and so
$H'$ has a matching $M'$ of size $\lfloor |H'|/3\rfloor$. One can then show that $M'$ contains at least $d$ edges from $H$, as desired. (We thank Peter Allen for
suggesting this trick.) 

However, the proof of Theorem~\ref{generalmatch} is only slightly simpler in the case when $d=\lfloor n/3 \rfloor$ (we do not need 
Claims~\ref{claimV'V0V0}--\ref{claimedges} in this case) and to show that the above trick works, one requires some extra calculations.%
\COMMENT{We are adding $a=(n-3d)/2-x/2$ vertices to $H$ where $x=0,1$. So now have $H'$ where
$\delta _1 (H')>\binom{n-1}{2}-\binom{n-d}{2}+\binom{a}{2}+a(n-1)=\binom{n+a-1}{2}-\binom{n-d}{2}=\binom{n+a-1}{2}-
\binom{n+a-\lfloor(n+a)/3 \rfloor}{2}.$ (The latter equation follows since 
$ n+a-\lfloor(n+a)/3 \rfloor = \lceil 2(n+a)/3 \rceil =\lceil n-d-x/3 \rceil=n-d$.) So $H'$ has a $\lfloor (n+a)/3 \rfloor$-matching.
This induces a matching in $H$ which covers all but at most $2a +(n+a \text{ mod } 3)$ vertices. 
Now $2a +(n+a \text{ mod } 3)=n-3d-x +(n+a \text{ mod } 3)$. If $n+a =3(n-d)/2 -x/2= 2 \text{ mod } 3$ then $-x =4 \text{ mod } 3$, i.e. $x \not = 
0,1$, a contradiction.  If $n+a =3(n-d)/2 -x/2= 0 \text{ mod } 3$ then $-x =0 \text{ mod } 3$, i.e. $x=0$. So $2a +(n+a \text{ mod } 3= n-3d$.
If $n+a =3(n-d)/2 -x/2= 1 \text{ mod } 3)$ then $-x =2 \text{ mod } 3$, i.e. $x=1$. But then $2a +(n+a \text{ mod } 3)= n-3d$. So in any case, at least
$3d$ vertices covered by matching, i.e. have a $d$-matching.
}

\medskip 

{\footnotesize \obeylines \parindent=0pt

\begin{tabular}{lll}
Daniela K\"{u}hn, Deryk Osthus &\ &  Andrew Treglown \\
School of Mathematics &\ & School of Mathematical Sciences \\
University of Birmingham &\ & Queen Mary, University of London \\ 
Edgbaston &\ & Mile End Road \\
Birmingham &\ & London\\
B15 2TT &\ & E1 4NS\\
UK &\ & UK
\end{tabular}
}
\begin{flushleft}
{\it{E-mail addresses}:\\
\tt{\{kuehn,osthus\}@maths.bham.ac.uk}, \tt{treglown@maths.qmul.ac.uk}}
\end{flushleft}

\end{document}